\documentclass[twoside,11pt]{amsart}

\usepackage{a4,amsthm,latexsym,amsmath,amssymb,amsfonts,euscript}

\newtheorem{theorem}{Theorem}    
\newtheorem{proposition}[theorem]{Proposition}  
  
\newtheorem{lemma}[theorem]{Lemma}

\def\R{\mathbb R}
\def\N{\mathbb N}
\def\Z{\mathbb Z}
\def\Q{\mathbb Q}

\def\Ret{{\EuScript R}}

\def\Re{{\mathcal R}}

\def\W{{\mathcal W}}

\newcommand{\eqdef}{\stackrel{\scriptscriptstyle\rm def}{=}}

\begin{document}

\title
{Local rates of Poincar{\'e} recurrence for rotations and weak mixing}

\author{J.-R. Chazottes}
\address{
Centre de Physique Th{\'e}orique,
CNRS-UMR 7644, 
{\'E}cole Polytechnique,
91128 Palaiseau Cedex, France.
}
\email{jeanrene@cpht.polytechnique.fr} 

\author{F. Durand}  
\address{Laboratoire Ami{\'e}nois
de Math{\'e}matiques Fondamentales  et
Appliqu{\'e}es, CNRS-UMR 6140, Universit{\'e} de Picardie
Jules Verne, 33 rue Saint Leu, 80039 Amiens Cedex 1, France.}
\email{fdurand@u-picardie.fr}

\begin{abstract}
We study the lower and upper local rates of Poincar{\'e} recurrence
of rotations on the circle by means of symbolic dynamics. 
As a consequence, we show that 
if the lower rate of Poincar{\'e} recurrence of an ergodic dynamical
system $(X,\EuScript{F} , \mu, T)$ is greater or equal to 1 $\mu$-almost everywhere, then
it is weakly mixing.
\end{abstract}


\maketitle

\section{Introduction and main result}\label{intro}

Let $T$, acting on the Lebesgue probability space $(X,\EuScript{F},\mu)$, be a
(non necessarily invertible) measure-preserving ergodic map.
Let $U$ be a measurable subset of $X$ and $x\in U$. Define the first return time of $x$ to $U$ by
$\tau_{U}(x)\eqdef\inf\{k\geq 1:T^{k}x\in U\}$. If $\mu(U)>0$ then $\tau_U(x)$ is finite
for $\mu$-almost every $x$ by Poincar{\'e} recurrence theorem.
Now define the  Poincar{\'e} recurrence of the set $U$ by 
$$
\tau(U)\eqdef \inf\{\tau_{U}(x):x\in U\}\;.
$$
It is easy to check that 
$$
\tau(U)=\inf\{k>0: T^{k}U\cap U\neq \emptyset \}=
\inf\{k>0: T^{-k}U\cap U\neq \emptyset \}\,.
$$
The notion of Poincar{\'e} recurrence of a set is used in
\cite{valya} to define a notion of dimension similar to Hausdorff
dimension where diameters of sets are replaced by their Poincar{\'e}
recurrence (see also 
\cite{acs,asuu,bruin,fuu,KPV,KM,psv}).

Suppose now that $\zeta$ is a finite measurable partition of $X$ and
denote, as usually, by $\zeta_n$ the partition
$\zeta\vee T^{-1}\zeta\vee\cdots\vee T^{-n+1}\zeta$
($n\geq 1$, $\zeta_1\eqdef\zeta$) and by $\zeta_n(x)$
the atom of this partition containing point $x$.
We define the (lower and upper) local rate of Poincar{\'e} recurrence
for the partition $\zeta$ respectively as follows:

$$
\underline{\Ret}_{\zeta}(x)\eqdef \liminf_{n\to\infty}
\frac{\tau(\zeta_{n}(x))}{n}\, ,\quad
\overline{\Ret}_{\zeta}(x)\eqdef \limsup_{n\to\infty}
\frac{\tau(\zeta_{n}(x))}{n}\, .
$$

(Of course, these quantities depend on the map $T$ but we omit this
dependence in the notation.)
If $x\in X$ is a periodic point, then obviously $\underline{\Ret}_{\zeta}(x)=\overline{\Ret}_{\zeta}(x)=0$.
A useful fact to be used later is that
both $\underline{\Ret}_{\zeta}$ and $\overline{\Ret}_{\zeta}$ are sub-invariant
functions, namely $\underline{\Ret}_{\zeta}\circ T\leq\underline{\Ret}_{\zeta}$ and
$\overline{\Ret}_{\zeta}\circ T\leq\overline{\Ret}_{\zeta}$.
This is because 
for any $n\geq 1$ and any $x\in X$, $\tau(\xi_{n-1}(Tx))\leq\tau(\xi_{n}(x))$.
If we assume $\mu$ is an ergodic probability measure, then it follows
by basic arguments that
$\underline{\Ret}_{\zeta}$ and $\overline{\Ret}_{\zeta}$ are $\mu$-a.e. constant
functions (see \cite{acs} for details). 

The definition of the (lower and upper) local rate of Poincar{\'e}
recurrence first appeared in \cite{hsv} in connection with the error
term in the approximation by the exponential law of the distribution
of rescaled return times to cylinder sets.
In \cite{acs-era} the authors prove that the following result, where
$h_\mu(T,\zeta)$ denotes the measure-theoretic entropy of
$(X,\EuScript{F},\mu,T)$ with respect to the partition $\zeta$.

\begin{theorem}
\label{Rinf->h>0}
Let $(X,\EuScript{F},\mu,T)$ be  as above and $\zeta$ a finite partition of $X$.
If $h_{\mu}(T,\zeta)>0$ then $\underline{\Ret}_{\zeta}(x)\geq 1$ $\mu$-almost everywhere.
\end{theorem}

It is equivalent to state that if $\underline{\Ret}_{\zeta}(x)< 1$
$\mu$-almost everywhere, then $h_{\mu}(T,\zeta)=0$.
In the case when $(X,T)$ has the specification property (think of a
topologically mixing subshift of finite type over a finite alphabet as
a typical example), one has
$
\overline{\Ret}_{\zeta}(x)\leq 1
$
for all $x\in X$, where $\zeta$ is the canonical partition labelled by
the alphabet (see \cite{acs}). Therefore, in this case $\underline{\Ret}_\zeta(x)=\overline{\Ret}_\zeta(x)=1$
$\mu$-almost everywhere where $\mu$ is any ergodic probability measure
such that $h_{\mu}(T,\zeta)>0$.
Outside such sets of full measure,
the specification property allows one to construct some points $x$ such that
$0<\underline{\Ret}_\zeta(x)=\overline{\Ret}_\zeta(x)<1$.
On another hand,
the positiveness of entropy is an unavoidable assumption in Theorem
\ref{Rinf->h>0}.
Indeed, 
it is shown in \cite{acs} that $\underline{\Ret}_{\zeta}(x)=0$ for Lebesgue
almost every $x$ for a special class of rotations where $\zeta$ is the canonical partition into two
atoms given by the rotation angle.
A natural question is thus to know
whether or not the converse to Theorem \ref{Rinf->h>0} holds.
In the course of the present note, we will provide an example,
namely the Morse system, showing that this is not the case.

The motivation of the present work was to figure out which 
property of an ergodic dynamical system $(X,\EuScript{F},\mu,T)$,
equipped with a partition $\zeta$, the property 
``$\underline{\Ret}_\zeta(x)\geq 1$ $\mu$-almost surely'' is related to. 
In this direction, we have the following result.

\begin{theorem}\label{main}
Let $T $, acting on the Lebesgue probability space $(X,\EuScript{F},\mu , T)$, be a measure-preserving ergodic map. If 
$\underline{\Ret}_{\zeta}(x)\geq 1$ for $\mu$-almost every $x\in X$ and
every non-trivial measurable partition $\zeta$ then $(X,\EuScript{F},\mu , T)$ is weakly
mixing.
\end{theorem}
(By `non-trivial' we mean that no atom of $\zeta$ has measure $0$ or $1$.)
For definitions and properties of the classical notions of mixing in
ergodic theory, we refer the reader to e.g. \cite{petersen}.

Our approach to prove this theorem is as follows. The point is that
a non weakly mixing system has a non trivial eigenvalue, hence a
measure-theoretical factor which is a rotation. Moreover, 
every measurable partition of the factor induces a
measurable partition of the original system. Therefore we are led
to show that the lower local rate of Poincar{\'e} recurrence for rotations
is strictly less than 1 for some partition. But the rotations are measure
theoretically isomorphic to Sturmian subshifts. Hence, it suffices to
prove that the lower rate of
Poincar{\'e} recurrence for Sturmian subshifts is strictly less than 1
for some partitions. In fact we prove:

\begin{theorem}\label{bounded-case}
Let $(\Omega_{\alpha} , S)$ be the Sturmian subshift generated by
 $\alpha \in \R \setminus \Q$ and  $\mu$ be its unique ergodic measure.
Let $\zeta$ be the partition $\{ [0] , [1] \}$.
\begin{enumerate}
\item
\label{equivalence}
The following
 statement are equivalent.
\begin{enumerate}
\item
The coefficients of the continued fraction of
$\alpha $ are bounded;
\item
$\underline{\Ret}_{\zeta }(x) > 0$ for $\mu$-almost every $x$;
\item
$\overline{\Ret}_{\zeta }(x) < \infty$ for $\mu$-almost every $x$.
\end{enumerate}
\item
\label{inferieura1}
Moreover, if the coefficients of the continued fraction of $\alpha $ are
 bounded, then $\underline{\Ret}_{\zeta }(x)< 1$ and $\overline{\Ret}_{\zeta }(x) > 1$
for $\mu$-almost every $x$.
\item
\label{rotation}
The same results hold for $([0,1[ , x\mapsto x+ \alpha \mod 1)$ and
     $\zeta = \{ [0,1-\alpha [ , [ 1-\alpha , 1[ \}$.
\end{enumerate}
\end{theorem}

After completing our work, we discovered the preprint
\cite{Ku} in which the author computes the precise values of $\underline{\Ret}$ and
$\overline{\Ret}$ for rotations. See also \cite{K} for substitutive subshifts.

\section{Local rates of recurrence for some zero entropy systems}\label{0-entropy}

In this section we establish some results about the local rate of recurrence for some zero entropy
systems that will be useful in the proof of Theorem \ref{main}.

Let $A$ be a finite alphabet. We endow $A^\Z$ and $A^\N$ with the product topology. Let $S:A^\Z \to A^\Z$ be the shift map: 
$(S x)_n = x_{n+1}$ where $x=(x_n)\in A^\Z$. If $X$ is a closed and $S$-invariant set of $A^\Z$ then $(X,S)$ is
called a subshift. It is a continuous map. For all word $u$ on the alphabet $A$, we call 
cylinder generated by $u$ the set $[u]=\{(x_n)\in X : x_0 x_1 \cdots x_{|u|-1} = u \}$.

\subsection{Linearly recurrent subshifts.}

Let $x=(x_n)$ be a sequence in $A^\Z$ where $A$ is a finite alphabet. We call
$L(x)$ the set of all finite words appearing in $x$. The length $|u|$ of a
word $u\in L(x)$ is the number of symbols in $u$. 
Let $u,w\in L(x)$.
For a subshift $(X,S)$ we set $L(X)=\bigcup_{x\in X} L(x)$. We say an
element $x$ of $A^\N$ or $A^\Z$ generates $(X,S)$ if $L(x) = L(X)$.
We say that $w$ is a return word
to $u$ of $x$ if
$wu$ belongs to $L (x)$, $u$ is a prefix of $wu$ and $u$ has exactly two occurrences in $wu$.
(There is a more general definition of return words in \cite{DHS}.) 
We say that $x$ is uniformly recurrent if all $u\in L(x)$ appears infinitely many times in $x$, and,
for all $u\in L(x)$ there exists $K_u$ such that for all return word $w$
to $u$, $|w|\leq K_u |u|$. 
Let us recall that if the subshift $(X,S)$ is minimal then all its points are uniformly recurrent.

We say that $x$ is {\rm linearly recurrent} (LR) (with constant $K\in\R\backslash\{0\}$) if it is uniformly
recurrent and if for all $u\in L(x)$
and all return word $w$ to $u$, we have $|w| \leq K|u|$.
We say that a subshift $(X,S)$ is LR (with constant $K$) if it is minimal and contains a LR sequence (with constant $K$). 
Notice that a minimal subshift is LR if and only if all its elements
are linearly recurrent.

Let $u$ be a word and $\alpha \in \R_+$. The prefix of length 
$\lfloor |u| \alpha \rfloor $ of the sequence $uuu \dots $ is denoted by
$u^\alpha$, where $\lfloor . \rfloor$ is the integer part map.

\begin{proposition}[\cite{DHS}]
\label{proplr}
Let $(X,S)$ be an aperiodic LR subshift with constant $K$. Then               
$X$ is $(K+1)$-power free
(i.e. $u^{K+1}\in  L(X)$ if and only if
$u=\emptyset$)
and
for all $ u\in  L(X)$ and for all $w \in\Re_{u}$ we have $(1/K)|u| < |w| $.   
\end{proposition}

Let $(X,S)$ be a subshift on the alphabet $A$. 
Let $\zeta $ be the partition into $1$-cylinders, i.e. $\zeta = \{ [a] ;
a \in A \}$.
There is a bijective correspondence between the
partition $\zeta_n$ and the set of the words of length $n$ of
$X$. Moreover, for all $x=(x_n)\in X$ and $n\in \N$ we have
$$
\tau (\zeta_n (x))
= 
\min \{ |w| : w\,\textup{is a return word to}\, x_0 x_1 \dots x_{n-1}   \}\, .
$$
Hence by Proposition \ref{proplr} we have the following result.

\begin{proposition}\label{ratelinrec}
Let $(X,S)$ be an aperiodic LR subshift with constant $K$ and $\zeta $ the
partition into $1$-cylinders.
Then for all $x\in X$:  
$$
\frac{1}{K}
\leq
\underline{\Ret}_{\zeta} (x) 
\leq
\overline{\Ret}_{\zeta} (x) 
\leq
K .
$$
\end{proposition}


Let $(X,S)$ be a subshift and $\zeta $ be the partition into $1$-cylinders.
If a subshift is $(K+1)$-power free then
$1/K \leq \underline{\Ret}_{\zeta} (x)$ for all $x\in X$. Hence if we define 
$$
\delta =  \inf \{ K\in\R\backslash\{0\}:(X,S) \hbox{ is } (K+1)-\hbox{power free} \}
$$
then it comes that  $1 / \delta \leq \underline{\Ret}_{\zeta} (x)$ for
all $x\in X$.
It is known \cite{Th} that $\delta = 1$ for the Morse sequence
(this is the fixed point of the substitution $\sigma$ defined by $\sigma (0) = 01$ and $\sigma (1) = 10$).
For the subshift $(X,S)$ of $\{0,1\}^{\N}$ generated by the $S$-orbit closure of the Morse sequence, one
has $\underline{\Ret}_{\zeta} (x)\geq 1$ for all $x\in X$, where $\zeta=\{[0],[1]\}$.
This gives an example of a zero-entropy dynamical system having a local rate
of Poincar{\'e} recurrence greater or equal to one. Therefore the converse to
Theorem \ref{Rinf->h>0} is false.

\subsection{Sturmian subshifts}

\label{sturmsubshift}
Let \( 0<\alpha <1 \) be an 
irrational number. We define the map 
\( R_{\alpha }:\left[ 0,1\right[ \rightarrow \left[ 0,1\right[  \) 
by 
$R_{\alpha }\left( t\right) =t+\alpha$
(mod 1) and the map 
$I_{\alpha }:\left[ 0,1\right[ \rightarrow \left\{ 0,1\right\}$
by 
$I_{\alpha }\left( t\right) =0 $
if 
$t\in \left[ 0,1-\alpha \right[  $
and 
$ I_{\alpha }\left( t\right) =1 $
otherwise. 
Let $\Omega _{\alpha }$ be the closure of the set $\left\{ \left. \left( I_{\alpha }
\left( R^{n}_{\alpha }\left( t\right) \right) \right) _{n\in \Z }\textrm{ }\right| \textrm{ }t\in
\left[ 0,1\right[ \textrm{ }\right\}$. 
The subshift 
$ \left( \Omega _{\alpha },S \right) $
is the {\it Sturmian subshift} generated by $\alpha$ and its elements are called {\it Sturmian 
sequences}.
There exists a factor map (see \cite{HM}) 
$\phi :\left( \Omega _{\alpha },S \right) \rightarrow \left( \left[ 0,1\right[ ,R_{\alpha }\right) $
such that 
$
\left| \phi ^{-1}\left( \left\{ \beta \right\} \right) \right| =2$
if
$\beta \in \left\{ \left. n\alpha \textrm{ }\right| \textrm{ }n\in \Z \textrm{ }\right\}$ 
and
$\left| \phi ^{-1}\left( \left\{ \beta \right\} \right) \right| =1$
otherwise.
Consequently $\phi$ is a measure-theoretic isomorphism.
It is well-known that $\left( \Omega _{\alpha },S \right) $
is a non-periodic uniquely ergodic minimal subshift. 

\begin{proposition}
[\cite{Du2,Du3}]
\label{lrsturm}
A Sturmian subshift $( \Omega _{\alpha },S)  $ is LR
if and only if
the coefficients of the continued fraction expansion of $\alpha$ are bounded.
\end{proposition}

In the sequel we will make use of the following
morphisms 
\( \rho _{n} \) 
and 
\( \gamma _{n} \), 
\( n\in \N \setminus \{0\} \), 
from 
\( \left\{ 0,1\right\}  \) 
to 
\( \left\{ 0,1\right\} ^{*} \) 
defined 
by
\[
\begin{array}{l}
\rho _{n}\left( 0\right) =01^{n+1}\\
\rho _{n}\left( 1\right) =01^{n}
\end{array}\textrm{ and }\begin{array}{l}
\gamma _{n}\left( 0\right) =10^{n+1}\\
\gamma _{n}\left( 1\right) =10^{n}
\end{array}\, .\]

\begin{proposition}\label{kappa}
\label{sturmdb} 
Let 
\( \left( \Omega_{\alpha},S\right)  \) 
be a Sturmian 
subshift. 
There exists a sequence
\( \left( \kappa_{n}\right) _{n\in \N } \) %
taking values in 
\( \left\{ \rho _{1},\gamma _{1},\rho _{2},\gamma _{2},\ldots \right\}  \) 
such that 
\begin{enumerate}
\item 
\( \displaystyle y=\lim _{n\rightarrow +\infty }\kappa _{1}\cdots \kappa _{n}\left( 00\cdots \right)  \) 
exists and generates 
\( \left( \Omega_{\alpha},S \right)  \);
\item
$ \left( \Omega_{\alpha},S\right)   $ is uniquely ergodic;
\item
$1 \leq \frac{| \kappa _{1}\cdots \kappa _{n} (0)|}{ |\kappa _{1}\cdots \kappa _{n} (1)| } \leq \frac{3}{2}$;

\item 
Let \( P_{0}=\{[0],[1]\}, \) and for \( n\geq 1, \) let 
$$
P_{n}=\left\{ \left. S ^{k}\kappa _{1}\cdots \kappa _{n}\left( \left[
        a\right] \right) \textrm{ }\right|
\textrm{ }0\leq k<\left| \kappa _{1}\cdots \kappa _{n}\left( a\right) \right| ,
\textrm{ }a\in \left\{ 0,1\right\} \textrm{ }\right\}
$$
is a partition of \( \Omega_{\alpha} \) with the following properties: 

\begin{enumerate}
\item 
\( \kappa_{1}\cdots \kappa_{n+1}\left( \left[ 0\right] \right) \cup \kappa
_{1}\cdots \kappa _{n+1}\left( \left[ 1\right] \right)
\subseteq \kappa _{1}\cdots \kappa _{n}\left( \left[ 0\right] \right)
\cup \kappa _{1}\cdots \kappa _{n}\left( \left[ 1\right] \right)  \), 
\item \( P_{n}\prec P_{n+1} \) 
as partitions,

\end{enumerate}
\end{enumerate}
\end{proposition}

\begin{proof}
The statements 1., 2.  and 4. are mainly due to Hedlund and Morse \cite{HM} (see \cite{DDM}). Point 3. follows from the definition
of $\rho_n$ and $\gamma_n$ given above.
\end{proof}


\section{Proof of Theorem \ref{main} and Theorem \ref{bounded-case}}


\subsection{Proof of Theorem \ref{bounded-case}}

In this subsection $\alpha \in [0,1[$ is an irrational number,
$(\Omega_\alpha , S)$ the Sturmian subshift it defines, $\mu$ its unique
ergodic measure and
$(\kappa_n)_{n\in \N}$ the sequence given by Proposition \ref{sturmdb}.

For all $n\in \N$ and all $a \in \{ 0,1 \}$ we set 
\begin{align}
\label{mu}
\mu_n (a) =  \sum_{0\leq k<\left| \kappa _{1}\cdots \kappa _{n}(a) \right| }
\mu \left( S^{k}\kappa _{1}\cdots \kappa _{n} ([a]) \right)
=
\left| \kappa _{1}\cdots \kappa _{n}(a) \right|
\mu \left( \kappa _{1}\cdots \kappa _{n} ([a]) \right) .
\end{align}
We remark that $\mu_n(0) + \mu_n (1) = 1$.



\begin{lemma}
\label{sturmrate} 
Let $n\in \N$ and $m\in \R$ with
$1\leq m  \leq |\kappa_n (0)| -2$. We set $\kappa_n (0) = ab^{|\kappa_n (0)| - 1}$. Let $U_n$ be the union of the sets
$S^k \kappa_{1}\cdots \kappa_{n} ([0])$, where

\begin{equation}
\label{ineq1}
|\kappa _{1}\cdots \kappa_{n-1} (a)| \leq k \leq | \kappa _{1}\cdots \kappa _{n} (0)| -\lfloor | \kappa _{1}\cdots \kappa_{n-1} (b)| (m+1) \rfloor ,
\end{equation}

and $V_n$ be the union of the sets
$S^l \kappa_{1}\cdots \kappa_{n} ([1])$, where

\begin{equation}
\label{ineq2}
|\kappa _{1}\cdots \kappa_{n-1} (a)| \leq l \leq | \kappa _{1}\cdots \kappa _{n} (1)| -\lfloor | \kappa _{1}\cdots \kappa_{n-1} (b)| m \rfloor .
\end{equation}

If $x\in U_n$ then there exists a word $v$ such that 
$x\in [v^{m+1}]$. If $x\in  V_n$ then there exists a word $v$ such that 
$
x\in [v^m] .
$
In both case $|v| = |\kappa_1 \cdots \kappa_{n-1} (b)|$.
Moreover for all $n\in \N$ 
$$
\mu (U_n ) \geq \mu_n(0) \left(\frac{|\kappa_n (0)| - m - 2}{3/2 + |\kappa_n (0)| - 1  } \right) \hbox{ and } \mu (V_n ) \geq \mu_n (1) \left(  \frac{ |\kappa_n (0)| -  m -2}{3/2 + |\kappa_n (0)| -2  } \right) .
$$
\end{lemma}

\begin{proof}
Let $n\in \N$ and $m\in [1 , |\kappa_n (0)| -2 ]$. Let $x\in S^k
\kappa_{1}\cdots \kappa_{n} [0] $ where $k\in
\N$ satisfies ($\ref{ineq1}$).
Let $w$ be the prefix of length $\lfloor |\kappa_0 \cdots \kappa_{n-1}
(b)|(m+1) \rfloor $ of $x$.
We have $x\in [w]$ and 
$$
\lfloor |\kappa_0 \cdots \kappa_{n-1}
(b)|(m+1) \rfloor 
=
|\kappa_0 \cdots \kappa_{n-1}
(b)|(\lfloor m \rfloor +1) 
+
\lfloor |\kappa_0 \cdots \kappa_{n-1}
(b)|(m - \lfloor m \rfloor ) \rfloor .
$$

The sequence $x$ belongs to
$S^k[ \kappa _{1}\cdots \kappa_{n-1} (a b^{|\kappa_n (0)| -1})]$. Consequently, the hypotheses on $k$ imply
that there exist two words $s$ and $p$ such that
$$
w = s (\kappa_0 \cdots \kappa_{n-1} (b))^{\lfloor m \rfloor } p (sp)^{ m - \lfloor m \rfloor }
=
(sp)^{m+1},
$$ 
where
$ps = \kappa_0 \cdots\kappa_{n-1} (b)$.
The other case can be treated in the same way.

From Proposition \ref{sturmdb} we deduce that
\begin{align*}
\mu (U_n) &
=
\mu_n (0) \frac{|\kappa _{1}\cdots \kappa _{n}(0)| - |\kappa _{1}\cdots \kappa _{n-1}(a)| - \lfloor |\kappa _{1}\cdots
\kappa _{n-1}(b)|(m+1) \rfloor +1  }{|\kappa _{1}\cdots \kappa _{n}(0)|} \\
& 
\geq
\mu_n (0)
\left(
\frac{
|\kappa _{1}\cdots \kappa _{n-1}(b)|(|\kappa_n (0)| - m - 2)
}
{
|\kappa _{1}\cdots \kappa _{n-1}(a)|+  |\kappa _{1}\cdots \kappa _{n-1}(b)| 
\left(
|\kappa_n (0)|-1
\right)
}
\right) \\
& \geq
\mu_n (0) \left(\frac{|\kappa_n (0)| - m - 2}{3/2 + |\kappa_n (0)| - 1
 } \right) .
\end{align*}
The same computations can be done for $V_n$.
\end{proof}

\begin{proof}[Proof of the statement \eqref{equivalence} of Theorem \ref{bounded-case}]
To prove that (a) is equivalent to (b) it suffices to prove that if the coefficients of the continued fraction of $\alpha $ are not bounded then 
$\underline{\Ret}_{\zeta } (x) = 0$ for $\mu$-almost every $x$. The
 other part of the proof follows from Proposition \ref{ratelinrec} and
 \ref{lrsturm}. Hence, $\underline{\Ret}_\zeta$ being $S$-invariant and
 $\mu$ ergodic, it is enough to prove that $\mu \left( \{ x;  \underline{\Ret}_{\zeta } (x) = 0 \} \right) > 0$.

Let $(U_n)$ and $(V_n)$ be the sequences of open sets given by Lemma \ref{sturmrate}. From
Proposition 1.1 in \cite{Du3} and Proposition \ref{lrsturm} there exists a strictly increasing sequence $(n_i)$ such that
$\lim_{i\rightarrow +\infty} |\kappa_{n_i} (0)| = +\infty $.
For all $i\in \N$ we set $m_i = \lfloor\sqrt{|\kappa_{n_i} (0)|
-2}\rfloor$.
For all $i\in \N$ and all $x\in U_{n_i} \cup V_{n_i}$, by Lemma \ref{sturmrate},
there exists $v_i$ such that $x$ belongs to the cylinder $[v_i^{m_i}]$ and
consequently 

$$
\frac{\tau \left(\zeta_{(m_i-1)|v_i|} (x)\right)}{(m_i-1)|v_i|} =
\frac{1}{m_i -1}  .
$$

Hence, if $x\in \cap_{j\in \N} \cup_{i\geq j} (U_{n_i}\cup V_{n_i} )$ then
$\underline{\Ret}_{\zeta} (x) = 0$. But, from Lemma \ref{sturmrate}, we
also have $\mu (\cap_{j\in \N} \cup_{i\geq j}(U_{n_i} \cup V_{n_i} )  ) \geq 2/3$. Thus, (a) is equivalent to (b).

\medskip

To prove that (a) is equivalent to (c) it suffices to prove that if the coefficients of the continued fraction of $\alpha $ are not bounded then 
$\overline{\Ret}_{\zeta } (x) = \infty$ for $\mu$-almost every $x$. The
 other part of the proof follows from Proposition \ref{ratelinrec} and
 \ref{lrsturm}. Hence, $\underline{\Ret}_\zeta$ being $S$-invariant and
 $\mu$ ergodic, it is enough to prove that $\mu \left( \{ x;  \overline{\Ret}_{\zeta } (x) \geq h  \} \right) > 0$ for all $h\geq 0$.

Let $h\geq 2$. Let $n\in \N$ be such that $\kappa_{n+1} (0)  =ab^{i+1}$
 with $i\geq 2$. We set $l_n = |\kappa _{1}\cdots \kappa_{n}(a)|$, $k_n = |\kappa _{1}\cdots \kappa_{n}(b^i)|/h$ and  $W_n = \cup_{1\leq k\leq k_n} 
S^{-k} \kappa _{1}\cdots \kappa_{n}([abb])$.

Take $x \in S^{-k} \kappa _{1}\cdots
\kappa_{n}([abb])$ with $1\leq k \leq k_n$. We remark $\kappa _{1}\cdots
\kappa_{n}([abb])$ is contained in $ \kappa _{1}\cdots
\kappa_{n+1}([a]) \cup \kappa _{1}\cdots
\kappa_{n+1}([b])$ (the proof is left to the reader). The words $\kappa_{n+1} (a)$
and $\kappa_{n+1} (b)$ end with the word $b^i$. Thus, we can write
$\kappa_1 \cdots \kappa_n (b) = uv$ in such a way that $v (\kappa_1
\cdots \kappa_n (b))^p \kappa_1
\cdots \kappa_n (ab^j) $ is a prefix of $x$ for some $p \leq i$, where $j=\max (i ,2 )$. It can be seen that
$\tau \left( \zeta_{k_n + l_n-1} (x) \right)$ is greater than $| \kappa_1 \cdots \kappa_n
(ab^j)|$. Consequently
\begin{align}
\label{retinfini}
\frac{\zeta_{k_n + l_n-1} (x)}{k_n + l_n }
\geq 
\frac{ j |\kappa_1 \cdots \kappa_n (b)| +  |\kappa_1 \cdots
\kappa_n (a)| }{|\kappa_1 \cdots \kappa_n (b^i)|/h +  |\kappa_1 \cdots
\kappa_n (a)|}
\geq \frac{j + 3/2}{i/h+ 3/2} =f(h,i,j) .
\end{align}
Hence if $x \in \W (h) = \cap_{j\in \N} \cup_{i\geq j} W_i$ then
$\overline{\Ret} (x) \geq f(h,i,j)$. But for all $n\in \N$ we have
\begin{align*}
\mu (W_n) & = k_n \mu (\kappa_1 \cdots \kappa_n ([abb])) 
\geq
k_n \mu \left( \kappa_1 \cdots \kappa_{n+1} ([0]) \cup \kappa_1 \cdots
 \kappa_{n+1} ([1]) \right) \\
 & =
\frac{i|\kappa_1 \cdots \kappa_{n} (b)|}{h}  
\left( 
\frac{\mu_{n+1}(0)}{|\kappa_1 \cdots \kappa_{n} (ab^{i+1}) |} 
+
\frac{\mu_{n+1}(1)}{|\kappa_1 \cdots \kappa_{n} (ab^{i}) |} 
\right)
\\
& 
\geq 
\frac{i}{h}
\frac{|\kappa_1 \cdots \kappa_{n} (b)|}{|\kappa_1 \cdots \kappa_{n}
 (ab^{i+1}) |} 
=
\frac{i}{h}
\frac{|\kappa_1 \cdots \kappa_{n} (b)|}{|\kappa_1 \cdots \kappa_{n} (a)|
 + (i+1) |\kappa_1 \cdots \kappa_{n} (b) |} 
\\
& 
\geq
\frac{i}{h}
\frac{1}{3/2 + (i+1)} 
\geq
\frac{2}{h}
\frac{1}{3/2 + (2+1)} 
\geq 
\frac{4}{9h} .
\end{align*}

Hence $\mu (\W (h)) \geq 4/9h 1 $. From
Proposition 1.1 in \cite{Du3} and Proposition \ref{lrsturm} there exists a strictly increasing sequence $(n_i)$ such that
$\lim_{i\rightarrow +\infty} |\kappa_{n_i} (0)| = +\infty $. Thus, using
\eqref{retinfini} it comes that $\overline{\Ret} (x) \geq h$ for all
 $x\in \W (h)$.  This proves that (a) is equivalent (c).
\end{proof}

\begin{proof}[Proof of the statement \eqref{inferieura1} of Theorem \ref{bounded-case}]

Using \eqref{retinfini} we
conclude that $\overline{\Ret} (x) > \frac{7}{5}$ for $\mu$-almost every $x\in
X$. 

Now we prove the other part of the statement. It suffices to prove that
 for some $\theta < 1$ we have $\mu ( \{ x : \underline{\Ret}_\zeta(x) \leq \theta \}) > 0$.

From the hypotheses the sequence $(|\kappa_n (0)| ; n\in \N)$ is
 bounded by some constant $K$. We will need the following lemma.

\begin{lemma}
\label{muzero}
For all $a\in \{ 0,1 \}$ and all $n\in \N$ we have $\mu_n (a) \geq 2/(3K+1)$.
\end{lemma}

\begin{proof}
Suppose $\kappa_{n+1} (0) = ab^{i+1}$ then one can check we have
\begin{align*}
\mu ( \kappa_1 \dots \kappa_{n} ([a])) 
= &
\mu ( \kappa_1 \dots \kappa_{n+1} ([0])) 
+
\mu ( \kappa_1 \dots \kappa_{n+1} ([1])) \hbox{ and } \\
\mu ( \kappa_1 \dots \kappa_{n} ([b])) 
= &
(i+1)\mu ( \kappa_1 \dots \kappa_{n+1} ([0])) 
+
i\mu ( \kappa_1 \dots \kappa_{n+1} ([1])) . \\
\end{align*}
Consequently using \eqref{mu} we obtain 
$$
\frac{3}{2i}
\geq
\frac{\mu_{n} (a)}{\mu_{n} (b)} 
\geq 
\frac{2}{3}
\frac{\mu ( \kappa_1 \dots \kappa_{n} ([a]))}{ \mu ( \kappa_1 \dots \kappa_{n} ([b]))  } 
\geq
\frac{2}{3(i+1)} .
$$
We conclude using the facts that $i\geq K$ and $\mu_n (a) + \mu_n (b) = 1$.
\end{proof}

We consider several cases. Suppose there exists a strictly increasing sequence $(n_i)$
 such that $|\kappa_{n_i} (0)|\geq 4$ for all $i\in \N$. 
We set $m_i = 3$, $i\in \N$.
Let $(U_{n_i})_{i\in \N}$ and $(V_{n_i})_{i\in \N}$ be the sequences of open sets given by
 Lemma \ref{sturmrate} and associated to the sequence $(m_i)_{i\in \N}$.
For all $i\in \N$ and all $x\in U_{n_i}$, by Lemma \ref{sturmrate},
there exists $v_i$ such that $x$ belongs to the cylinder $[v_i^{5/2}]$
and $\lim_{i\to \infty} |v_i| = + \infty$. Consequently 

$$
\frac{\tau \left(\zeta_{(3/2)|v_i|-1} (x)\right)}{ \lfloor (3/2)|v_i| \rfloor } \leq
\frac{|v_i|}{(3/2)|v_i|} = \frac{2}{3} ,
$$

and, if $x\in \cap_{j\in \N} \cup_{i\geq j}U_{n_i}$ then
$\underline{\Ret}_{\zeta} (x) \leq 2/3$. Lemma \ref{muzero} and Lemma \ref{sturmrate} imply
$$
\mu (U_{n_i})
\geq
\mu_n (0) \frac{|\kappa_{n_i} (0)| - 7/2}{1/2 + |\kappa_{n_i} (0)| } 
\geq
\frac{2(K - 7/2)}{(3K+1)(1/2 + K)} > 0.
$$
Hence $\mu (\cap_{j\in \N} \cup_{i\geq j}U_{n_i}) > 0$ and  $\mu ( \{ x : \underline{\Ret}_\zeta(x) \leq 2/3\}) > 0$.

It remains to treat the following case: There exists $i_0$ such that for all $i\geq i_0$, $|\kappa_i (0)| = 3$.
For all $n\in \N$ we define $W_n$ to be the union of the sets
$S^k \kappa_{1}\cdots \kappa_{2n} [0]$, where

\begin{equation}
\label{ineq3}
|\kappa _{1}\cdots \kappa_{2n-2} (01)| \leq k \leq | \kappa _{1}\cdots \kappa_{2n-2} (01)| 
+\lfloor | \kappa _{1}\cdots \kappa_{2n-2} (1)|/2 \rfloor .
\end{equation}

Let $x\in W_n$. We consider four cases. 

First case: Suppose $\kappa_{2n-1} \kappa_{2n} = \rho_1^2$. We have 
$\kappa_{2n-1} \kappa_{2n} (0) = 0110101$. Hence $x$ belongs to $T^k [
\kappa_1 \kappa_2 \cdots \kappa_{2n-2} (0110101)]$ for some $k$
satisfying (\ref{ineq3}).
We can write $\kappa_1
\kappa_2 \cdots \kappa_{2n-2} (1) = uv$, with $|v| \geq \lfloor | \kappa _{1}\cdots \kappa_{2n-2} (1)|/2 \rfloor $, in such a way that the word
$$
v \ \kappa_1 \kappa_2 \cdots \kappa_{2n-2} (0) \ uv \ \kappa_1 \kappa_2
\cdots \kappa_{2n-2} (0) \  uv
$$ 
is a prefix of $x$.
We set $k_n = | \kappa _{1}\cdots \kappa_{2n-2} (01)| + |v|$. We have 

\begin{align*}
\frac{
\tau 
\left(
\zeta_{k_n-1} (x)
\right)
}{k_n} 
= & 
\frac{| \kappa _{1}\cdots \kappa_{2n-2} (01)| }{k_n} 
\leq
\frac{| \kappa _{1}\cdots \kappa_{2n-2} (01)| }{| \kappa _{1}\cdots \kappa_{2n-2} (01)| + 
| \kappa _{1}\cdots \kappa_{2n-2} (1)|/2 -1 }  \\
\leq & 
\theta_1 =
\frac{2}{2+1/5} < 1,
\end{align*}

Second case: $\kappa_{2n-1} \kappa_{2n} = \gamma_1 \rho_1$. We have 
$\kappa_{2n-1} \kappa_{2n} (0) = 1001010$. As previously we obtain that there exists $\theta_2 < 1$ such that 

$$
\frac{
\tau 
\left(
\zeta_{k_n-1} (x)
\right)
}
{k_n} 
\leq
\theta_2 < 1 .
$$

Third case: $\kappa_{2n-1} \kappa_{2n} = \gamma_1 ^2$. 
We have 
$\kappa_{2n-1} \kappa_{2n} (0) = 10100100$. 
We write $\kappa_1 \kappa_2 \cdots \kappa_{2n-2} (1) = uv$ such that 
$v\kappa_1 \kappa_2 \cdots \kappa_{2n-2} (00) uv \kappa_1 \kappa_2
\cdots \kappa_{2n-2} (00) $ is a prefix of $x$ and 
$|v| \geq \lfloor | \kappa _{1}\cdots \kappa_{2n-2} (1)|/2 \rfloor $.
But the images of 0 and 1 by $\kappa_{2n-1} \kappa_{2n} $ begin with the letter $1$. Furthermore the word
$$
v \ \kappa_1 \kappa_2 \cdots \kappa_{2n-2} (00) \ uv \ \kappa_1 \kappa_2 \cdots
\kappa_{2n-2} (00) \ uv
$$ 
is a prefix of $x$.
We set $l_n = | \kappa _{1}\cdots \kappa_{2n-2} (001)| + |v|$.
We have 
\begin{align*}
\frac{
\tau 
\left(
\zeta_{l_n-1} (x)
\right)}
{l_n}  = &
\frac{| \kappa _{1}\cdots \kappa_{2n-2} (001)| }{l_n}  \\
\leq &
\frac{| \kappa _{1}\cdots \kappa_{2n-2} (001)| }{| \kappa _{1}\cdots \kappa_{2n-2} (001)| + 
 | \kappa _{1}\cdots \kappa_{2n-2} (1)|/2 -1 }  
\leq \theta_3 = \frac{2}{2+1/8} < 1.
\end{align*}

Fourth case: $\kappa_{2n-1} \kappa_{2n} = \rho_1 \gamma_1$. We have 
$\kappa_{2n-1} \kappa_{2n} (0) = 01011011$. Proceeding as in the third case we obtain that there exists $\theta_4 < 1$ such that 

$$
\frac{
\tau 
\left(
\zeta_{m_n-1} (x)
\right)}
{m_n} 
\leq
\theta_4<1  ,
$$

where $m_n = | \kappa _{1}\cdots \kappa_{2n-2} (011)| + |v|$.

\medskip

To conclude we set $\theta = \max \{ \theta_1, \theta_2 , \theta_3 ,\theta_4 \}$ and
 we remark we have
\begin{align*}
\mu (W_n)
\geq & \mu_n (0) 
\frac{
|\kappa _{1}\cdots \kappa_{2n-2}(1)|
}
{
2|\kappa _{1}\cdots \kappa _{2n}(0)|
} \\
\geq &
\mu_n (0)
\frac{
|\kappa _{1}\cdots \kappa_{2n-2}(1)|}
{
5|\kappa _{1}\cdots \kappa _{2n-2}(0)| 
+ 5|\kappa _{1}\cdots \kappa _{2n-2}(1)|
}
\geq
\frac{4}{25(3K+1)}.
\end{align*}

Hence $\mu (\cap_{j\in \N} \cup_{i\geq j}W_i) > 0$ and  $\mu ( \{ x : \underline{\Ret}_\zeta(x) \leq \theta\}) > 0$.
\end{proof}

\begin{proof}[Proof of the statement \eqref{rotation}]
It suffices to remark that $[0] = \phi^{-1} ([0, 1-\alpha[ )$ and
$[1] = \phi^{-1} ([1-\alpha,1[ )$ (where $\phi$ is the
 measure-theoretical isomorphism given in Subsection \ref{sturmsubshift}). 
\end{proof}

\subsection{Proof of Theorem \protect{\ref{main}}}

If the system is not weakly mixing then it has a nontrivial
eigenvalue $\exp(2i\pi \alpha )$, therefore a measure-theoretical factor which is a
rotation: $([0,1[ , R_\alpha)$.
Furthermore, every measurable partition of the factor induces a
measurable partition of the original system. 
If $\alpha$ is a rational number then clearly there exists
a non-trivial partition $\zeta$ of $[0,1[$ so that $\underline{\Ret}_\zeta(x)=0$
for $\mu$-almost every $x$. 
If $\alpha$ is irrational then by Theorem \ref{bounded-case} there
is a non-trivial partition $\zeta$ of $[0,1 [$ such that 
$\underline{\Ret}_\zeta(x) = 0$ for $\mu$-almost every $x$ if the coefficients
of the continued fration of $\alpha$ are not bounded and 
$0<\underline{\Ret}_\zeta(x)< 1$ for $\mu$-almost every $x$ otherwise. This ends the
proof. \hfill $\qed$

\medskip

{\bf Acknowledgments.} 
We kindly acknowledge the CMM in Santiago,
Chile, for the support. We also thank V. Afraimovich for stimulating 
suggestions when we met him at the workshop on Dynamics and Randomness
held at Santiago in 2000.


\end{document}